\documentclass{article}
\usepackage[letterpaper,left=1.2in,right=0.8in,top=1in,bottom=1in]{geometry}
\usepackage{stmaryrd,hyperref,amsmath,amsthm,amssymb,url,soul,enumerate}
\hypersetup{pdfborder={0 0 0}}
\usepackage[capitalise, noabbrev, nameinlink]{cleveref}
\usepackage{graphicx,xcolor,tikz,comment}
\usetikzlibrary{positioning,arrows.meta,calc}
\newdimen\Spt
\newdimen\Lpt
\definecolor{hlt2Fill}{rgb}{.8,.8,.8}
\usetikzlibrary{shapes}

\tikzset{construction graph settings/.style={semithick,every node/.append style={node font=\footnotesize,draw,inner sep=2.5pt},vertex/.style={minimum size=10.5pt,circle},operation/.style={diamond,outer sep=-0.2pt,minimum size=21pt,inner sep=1pt},construction edges/.style={very thick}}}

\usepackage[mathlines]{lineno}

\usepackage{color}
\definecolor{red}{rgb}{1,0,0}
\def\red{\color{red}}

\definecolor{grn}{rgb}{0,.8,0}

\DeclareMathOperator{\Z}{Z}
\newcommand{\zbar}{\overline\Z}
\DeclareMathOperator{\Zp}{Z_+}
\newcommand{\zpbar}{\overline\Zp}
\DeclareMathOperator{\pt}{pt}
\DeclareMathOperator{\ptp}{pt_+}
\newcommand{\PTbar}{\overline{\pt}}
\newcommand{\PTpbar}{\overline{\pt_+}}
\newcommand{\PTul}{\underline{\pt}}
\newcommand{\PTpul}{\underline{\pt_+}}
\DeclareMathOperator{\PT}{PT}
\DeclareMathOperator{\PTLV}{PTS}
\DeclareMathOperator{\PTpLV}{PTS_+}
\DeclareMathOperator{\EPTLV}{EPTS}
\DeclareMathOperator{\EPTpLV}{EPTS_+}
\newcommand{\thz}{\operatorname{th}}

\def\GraphSix#1"{\begingroup(graph6string:\nolinebreak[3]\ \verb"\aftergroup\endgroup\aftergroup)}

\numberwithin{figure}{section}
\newcommand{\Gc}{\overline{G}}
\newcommand{\du}{\mathbin{\,\dot{\cup}\,}}
\newcommand{\noi}{\noindent}

\newcommand{\bit}{\begin{itemize}}
\newcommand{\eit}{\end{itemize}}
\newcommand{\ben}{\begin{enumerate}}
\newcommand{\een}{\end{enumerate}}
\newcommand{\beq}{\begin{equation}}
\newcommand{\eeq}{\end{equation}}
\newcommand{\bea}{\begin{eqnarray*}}
\newcommand{\eea}{\end{eqnarray*}}
\newcommand{\bean}{\begin{eqnarray}}
\newcommand{\eean}{\end{eqnarray}}
\newcommand{\bpf}{\begin{proof}}
\newcommand{\epf}{\end{proof}\ms}
\newcommand{\bmt}{\begin{bmatrix}}
\newcommand{\emt}{\end{bmatrix}}
\newcommand{\ms}{\medskip}

\newcommand{\beqa}{\begin{array}}
\newcommand{\eeqa}{\end{array}}

\newcommand{\lc}{\left\lceil}
\newcommand{\rc}{\right\rceil}
\newcommand{\lf}{\left\lfloor}
\newcommand{\rf}{\right\rfloor}

\def\sgap#1{{\mathsf S}_{#1}}

\newtheorem{theorem}{Theorem}[section]
\newtheorem{corollary}[theorem]{Corollary}
\newtheorem{lemma}[theorem]{Lemma}
\newtheorem{observation}[theorem]{Observation}
\newtheorem{proposition}[theorem]{Proposition}

\newtheorem{conjecture}[theorem]{Conjecture}

\theoremstyle{definition}
\newtheorem{definition}[theorem]{Definition}
\newtheorem{example}[theorem]{Example}

\definecolor{purp}{rgb}{.5,0,.5}

\definecolor{dgreen}{rgb}{0,0.6,0.4}

\title{Zero forcing propagation time intervals and graphs with fixed propagation time}
\author{Daniela Ferrero\thanks{Department of Mathematics, Texas State University, San Marcos, TX 78666, USA. (dferrero@txstate.edu).}\and H.~Tracy Hall\thanks{Provo,  UT 84606, USA (h.tracy@gmail.com).}\and Leslie Hogben\thanks{American Institute of Mathematics, Pasadena, CA 91125, USA (hogben@aimath.org); Department of Mathematics, Iowa State University, Ames, IA 50011, USA; Department of Mathematics, Purdue University, West Lafayette, IN 47907, USA.} \and Mark Hunnell\thanks{Department of Mathematics, Winston-Salem State University, Winston-Salem, NC 27110, USA.\break (hunnellm@wssu.edu).}
\and Ben Small\thanks{University Place, WA 98466, USA. (bentsm@gmail.com).}}

\begin{document}

\maketitle
%\linenumbers

%%%%%%%%%%%%%%%%%%%%%%%%%%%%%%%%
\begin{abstract} 
 Zero forcing in a graph refers to the evolution of vertex states under repeated application of a color change rule. Typically the states are chosen to be blue and white, and a forcing set is an initial set of blue vertices such that all of the vertices are blue at the end of the process.  In this context, the propagation time of a set in a graph is the number of iterations of the color change rule required to have all vertices blue, performing independent color changes simultaneously.  Different minimal forcing sets need not have the same propagation time, and we study the realizability of specific  integers as propagation times of minimal forcing sets in graphs for two of the most well-studied color change rules (standard and positive semidefinite). Particular attention is paid to the case where all minimal forcing sets have the same propagation time, and we term this phenomenon fixed propagation time. For each of the two variants, we present a general form of graphs all of which have fixed propagation time equal to one. We conjecture that these are the only such graphs and prove the conjectures for joins of graphs. Families of graphs with longer fixed propagation time for standard forcing are exhibited, and it is shown that such graphs do not exist for positive semidefinite forcing.

\end{abstract}

\noi {\bf Keywords} propagation time; zero forcing; PSD forcing; fixed propagation time

\noi{\bf AMS subject classification}  05C69,  05C50, 05C57

%%%%%%%%%%%%%%%%%%%%%%%%%%%%%%%%%%%%%%%%%%%

%==========================
%
% Section
%
%==========================

\section{Introduction}

 Zero forcing describes a dynamic coloring process on a graph, where initially each vertex is assigned blue or white.  White vertices can become blue by applying a color change rule, which is iterated until a final coloring of the graph is obtained.  The version of the color change rule determines the variant of zero forcing and here we discuss both the standard and positive semidefinite (PSD) variants of zero forcing.  A set $B$ is a \emph{(standard or PSD) forcing set} for a graph $G$ if, starting with exactly the vertices in $B$ blue, repeated applications of the (standard or PSD) color change rule can color all vertices of $G$ blue.   Precise definitions of the color change rules and other terms used are given later in this introduction.

Zero forcing arose about twenty years ago in numerous applications, including control of quantum systems, the Inverse Eigenvalue Problem of a Graph, and graph searching, and there has been extensive study of both standard and PSD zero forcing (see \cite{HLS22book} and references therein).  

Initially, the main focus of the study of zero forcing was to determine the size of a minimum (standard or PSD)  forcing set, called the \emph{standard forcing number} $\Z(G)$ or the \emph{PSD  forcing number} $\Zp(G)$. 
Subsequently, interest developed in the number of iterations (called \emph{rounds} or \emph{time-steps}) needed to color all vertices blue when performing all independently possible forces simultaneously, called the \emph{propagation time}. When starting with a specific forcing set $B\subseteq V(G)$, the propagation time of this set is denoted by $\pt(G,B)$ or $\ptp(G,B)$.  The \emph{standard propagation time  of $G$}, $\pt(G)$, is the minimum propagation time of $\pt(G,B)$ among all minimum standard forcing sets $B$, and analogously for $\ptp(G)$ \cite{proptime, PSDproptime}.

In the last five years there has been much more interest in the study of all forcing sets.  This includes determination of the maximum cardinality of a minimal forcing set, called the \emph{upper (standard or PSD) forcing number} and the reconfiguration of forcing sets.
Reconfiguration examines relationships among solutions to a problem, with solutions considered abstractly as the vertices of a \emph{reconfiguration graph}.
A reconfiguration rule is given for the types of small changes that are allowed when trying to transform one solution to another, %and each such permissible reconfiguration defines an edge
defining the edges of the reconfiguration graph.
Depending on the rule chosen, a given pair of solutions may or may not belong to the same connected component of the reconfiguration graph,
and the nature of these connected components, for various rules and for various types of zero forcing, is an active area of study.
An introduction to reconfiguration, in both its broader context and its particular application to zero forcing, is provided by the bibliography of
\cite{bong2025psdskew} for example, which studies certain reconfiguration graphs for PSD forcing and skew forcing.

Propagation time does not lend itself directly to a reconfiguration rule on the solution sets themselves, but an analogous question is raised if one considers a set of abstract vertices defined by the different possible propagation times, with an edge between times that differ  exactly by one.
Thus when there is a missing value between the maximum and minimum propagation times over minimum forcing sets, i.e., an intermediate value that cannot be realized as the propagation time of a minimum forcing set, this corresponds to a disconnected reconfiguration graph.  Similar questions can be asked for minimal forcing
sets.  Results related to these ideas are presented in Sections \ref{s:PSDgap} and \ref{s:bigfixedpt}.

The \emph{(standard or PSD) propagation time interval} is the interval  of integers between the minimum of the (standard or PSD) propagation times of minimum (standard or PSD) forcing sets of $G$ and the maximum of the  (standard or PSD) propagation times of minimum  (standard or PSD) forcing sets of $G$; the propagation time interval is \emph{full} if every integer in the  propagation time interval is achievable by some minimum forcing set  \cite{proptime, PSDproptime}. If a propagation time interval is full, then one can move between any possible propagation times of minimum forcing sets in increments of one (possibly requiring major changes to the minimum forcing sets).

\begin{definition}
The \emph{(standard or PSD) propagation time set} of a graph $G$, denoted by $\PTLV(G)$ or $\PTpLV(G)$,  is the set of (standard or PSD) propagation times realized by minimum (standard or PSD) forcing sets.  
 
 For a graph $G$, $\PTLV(G)$ has a  \emph{gap} if $k,k+r\in \PTLV(G)$ for $r\ge 2$ but $k+1\not\in \PTLV(G)$, and similarly for  $\PTpLV(G)$.     
\end{definition}

The  propagation time interval of $G$ is full if and only if the set of integers in the propagation time interval equals the   propagation time set of $G$ of the same type.  
It was shown in \cite{proptime}  that the standard propagation time interval need not be full.   In 2016 Warnberg conjectured that every  PSD propagation time interval is full \cite{PSDproptime}. In 2020 Ben Small constructed an example of a graph that does not have a full PSD propagation time interval (shown in Figure \ref{f:gapex}), and the existence of such a graph was announced in \cite{HLS22book}. In Section \ref{s:PSDgap} this example is presented and it is proved that
%theory behind the original example is presented, leading to an entire family of graphs that do not 
its PSD propagation time interval is not full.

Now that minimal  forcing sets are studied, it is of interest to study their propagation times, and to minimize and maximize those values, as in the next definition.

\begin{definition} For a graph $G$, the \emph{upper standard propagation time $\PTbar(G)$} (respectively, the \emph{upper PSD propagation time $\PTpbar(G)$}) is the maximum propagation time of any minimal standard (respectively, PSD) zero forcing set. The \emph{lower standard propagation time $\PTul(G)$} (respectively, the \emph{lower PSD propagation time $\PTpbar(G)$}) is the minimum propagation time of any minimal standard (respectively, PSD) zero forcing set. 

 The \emph{expanded standard propagation time interval}  of a graph $G$ is defined to be $[\PTul(G),\PTbar(G)]=\{\PTpul(G),\PTul(G)+1,\dots,\PTbar(G)-1,\PTbar(G)\} $  and the  \emph{expanded PSD propagation time interval} is defined analogously. An expanded propagation time interval is \emph{full} if every integer in the propagation time interval is realized by some minimal forcing set.  
\end{definition}

 We denote the maximum propagation time of a \emph{minimal} forcing set by $\PTbar(G)$ and $\PTbar_+(G)$. Note that in the study of propagation time intervals (of minimum forcing sets) in the literature, the notations $\PT(G)$ and $\PT_+(G)$ are used to denote the maximum propagation time of a \emph{minimum} forcing set.  However, we follow the  reconfiguration literature for zero forcing related parameters in putting a bar over the name of the main propagation  parameter (which is the minimum time of a minimum set) to denote the maximum propagation time of a \emph{minimal} forcing set. The minimum propagation time of a  minimal forcing set is denoted  by $\PTul(G)$  and $\PTul_+(G)$, respectively. 
The possibility that $\PTul(G) < \pt(G)$ (and similarly for PSD propagation time) %The distinction between propagation times for minimum and 
occurs because we use minimal rather than minimum forcing sets.  A path and standard forcing provides an example illustrating this: It is well-known that $\pt(P_n)=n-1$, but $\PTul(P_n)=\lf\frac{n-2}2\rf$ using a set of two adjacent vertices in the center of the path.

 \begin{definition} The \emph{expanded (standard or PSD) propagation time set} of a graph $G$, denoted by $\EPTLV(G)$ or $\EPTpLV(G)$,  is the set of (standard or PSD) propagation times realized by minimal (standard or PSD) forcing sets.  
 
 For a graph $G$, $\EPTLV(G)$ has a  \emph{gap} if $k,k+r\in \EPTLV(G)$ for $r\ge 2$ but $k+1\not\in \EPTLV(G)$, and similarly for  $\EPTpLV(G)$. 
\end{definition}

As is the case for the propagation time interval, the expanded propagation time interval of $G$ is full if and only if  the set of integers in the propagation time interval equals the  expanded  propagation time set of $G$ of the same type.  
 The graph in Figure \ref{f:gapex} does not have a full expanded PSD propagation time interval.
 A gap in an expanded  propagation time set identifies a place where we cannot change from one propagation time to  another with an increment of one using minimal forcing sets.

The simplest type of expanded propagation time interval, which is guaranteed to be full, is one where the upper  and lower propagation times are equal.

\begin{definition} A graph $G$ has \emph{fixed propagation time} (respectively, \emph{fixed PSD propagation time}) if $ \PTul(G)=\PTbar(G)$ (respectively, $ \PTpul(G)=\PTpbar(G)$).
\end{definition}

 Equivalently, a graph $G$ has fixed propagation time if and only if the expanded propagation time interval has exactly one element. Note that $|\EPTLV(G)|=1$ implies $|\PTLV(G)|=1$, but the reverse implication does not hold, as evidenced by a path on at least four vertices. 
As discussed in Section \ref{s:bigfixedpt}, it is easy to find examples of graphs that have fixed standard propagation times greater than one. In that section we also show that a fixed PSD propagation time greater than one is not possible.

Fixed propagation time equal to one, i.e., the study of graphs $G$ such that   $\PTul(G)=\PTbar(G)=1$ or $\PTpul(G)=\PTpbar(G)=1$, is of particular interest because it is equivalent to every forcing set $B\ne V(G)$ having propagation time  equal to one (not just the minimal forcing sets). It is also simpler, in part because $\PTbar(G)=1$ implies $\PTul(G)=1$, and similarly for PSD. 
In Section \ref{s:fastjoins} we show that graphs with two specific forms have fixed standard or PSD propagation time equal to one; we conjecture (see Conjectures \ref{c:psd} and \ref{c:std}) that every graph having fixed standard or PSD propagation time equal to one has the associated form.  %{\blu In Section \ref{s:joinisfast} we show that these conjectures are true for graphs that are joins.}

 Next we define these specific forms. Let $G_1$ and $G_2$ be two vertex disjoint graphs (i.e. $V(G_1)\cap V(G_2)=\emptyset$). The \emph{(disjoint) union} of $G_1$ and $G_2$ is the graph $G=G_1\du G_2$ with $V(G)=V(G_1)\cup V(G_2)$ and $E(G)=E(G_1)\cup E(G_2)$. The \emph{join} of disjoint $G_1$ and $G_2$ is the graph $G_1\vee G_2$ having 
$V(G)=V(G_1)\cup V(G_2)$ and $E(G)=E(G_1)\cup E(G_2)\cup \{v_1v_2:v_1\in V(G_1), v_2\in V(G_2)\}$.

\begin{definition}
  A graph $G$ of order $n\geq 2$ is a \emph{PSD fast join} if $G\cong K_n$  (i.e. $G\cong K_{n-t}\vee K_t$ for $1\leq t\leq n-1$) or $G\cong (K_{n_1}\du K_{m_1})\vee \ldots \vee (K_{n_t}\du K_{m_t})$ for $t\geq 2$ and positive integers $n_i$ and $m_i$ such that $n=\sum_{i=1}^tn_i+m_i$.

A graph $G$ of order $n\geq 2$ is a \emph{standard fast join} if $G\cong K_n$ or $G=\bigvee_{i=1}^t G_i$ where $t\geq 2$ and each graph $G_i$ satisfies $G_i\cong n_iK_1$ or $G_i\cong K_{n_i}\du K_{m_i}$ with $n_i,m_i\geq 2$. 
\end{definition}

In Section \ref{s:fastjoins} we establish that  fast joins have fixed propagation time equal to one:

\begin{theorem}
%\label{t:ConnPSDClassHasFPT1}
If $G$ is a PSD fast join, then $\PTpbar(G)=1$.  If $G$ is a standard fast join, then $\PTbar(G)=1$.
\end{theorem}

We conjecture  the converses of  the statements in this theorem  are  true %\lh{grammar problem as written} 
and in Section \ref{s:joinisfast} we prove these conjectures for graphs that are joins or have sufficiently large forcing numbers.

\begin{conjecture}\label{c:psd}
Let $G$ be a connected graph. % of order at least two. 
If   $\PTpbar(G)=1$, then  $G$ is a PSD fast join.
\end{conjecture}

\begin{conjecture}\label{c:std}
Let $G$ be a connected graph. % of order at least two. 
If   $\PTbar(G)=1$, then  $G$ is a standard fast join.
\end{conjecture}

We conclude this introduction with precise definitions of the notation and terminology.

A graph $G$ is a pair $(V(G),E(G))$, where $V(G)$ is the set of vertices and $E(G)$ is the set of edges.  The elements of $E(G)$ are  subsets of vertices of cardinality two; we denote the edge between distinct vertices $u$ and $v$ by $uv$ or $vu$. If  $uv \in E(G)$, then $u$ and $v$ are called \emph{neighbors}. The \emph{open neighborhood} of $u$ in $G$, denoted $N_G(u)$, is the set of all neighbors of $u$ in $G$, while the \emph{closed neighborhood} of $u$ in $G$ is $N_G[u]=N_G(u)\cup \{u\}$. The \emph{degree} of $u$  in $G$ is the cardinality of $N_G(u)$ and is denoted be $\deg_G(u)$. A vertex $u$ satisfying $\deg_G(u)=1$ is called a \emph{leaf}, while if $|V(G)|=n$ and $\deg_G(u) =n-1$, then $u$ is called a \emph{universal vertex}.

We  consider two color change rules, CCR-$\Z$ for standard zero forcing and CCR-$\Zp$ for PSD forcing, and use  CCR-$R$ to denote either one of these.   
\begin{itemize}
    \item CCR-$\Z$: If a blue vertex $u$ has a unique white neighbor $w$, then $u$ can force $w$ to become blue.
    \item CCR-$\Zp$: If $B$ is the set of currently blue vertices, $C$ is a component of $G-B$, and $u$ is a blue vertex such that $N_G(u) \cap V(C)=\{w\}$, %for some white vertex $v$
    then $u$ can force $w$ to become blue.
\end{itemize}

For $B\subseteq V(G)$ and a color change rule CCR-$R$, we define two sequences of sets: %, the set $B^{(i)}$ of vertices that turn blue in round $i$   and  the set $B^{[i]}$ of vertices that are blue after round $i$. That is, 
$B^{[0]}=B^{(0)}=B$ is the  set of initially blue vertices.  Assume $B^{(i)}$ and $B^{[i]}$ have been constructed. Then
 \[B^{(i+1)}=\{w:\mbox{$w$ can be forced by some $v$ (given  $B^{[i]}$  blue)}\}\ \ \mbox{  and }\ \ B^{[i+1]}=B^{[i]}\cup B^{(i+1)}.\]  
 The vertices in $B^{(i)}$ are said to be forced in \emph{round $i$}. Let $p$ denote the greatest integer such that  $B^{(p)}\ne \emptyset$ (note $B^{(i)}= \emptyset$ implies $B^{(i+1)}= \emptyset$). Then $B$ is a  \emph{forcing set} of $G$ if and only if  $B^{[p]}=V(G)$. 

 A \emph{minimum} forcing set $B$ is a forcing set of minimum cardinality, {i.e.}, if $|B'|<|B|$, then $B'$ is not a forcing set.  A \emph{minimal} forcing set $B$ is a forcing set such that no proper subset is a forcing set, i.e., if $B'\subsetneq B$, then $B'$ is not a forcing set.  
 The (standard and PSD) zero forcing numbers of $G$, denoted by  $\Z(G)$ and $\Zp(G)$, are the minimum cardinalities among (standard or PSD) forcing sets. Note that every minimum forcing set is a minimal forcing set, but not all minimal forcing sets are minimum forcing sets.  Every forcing set contains a subset (not necessarily proper) that is a minimal forcing set, but not necessarily a subset that is a minimum forcing set.
 
 When $B$ is a  forcing set and  $p$ is the greatest integer such that  $B^{(p)}\ne \emptyset$, then  $p$ is called the \emph{propagation time} of $B$ in $G$ and denoted $\pt(G,B)$ and $\ptp(G,B)$ for standard and PSD zero forcing sets, respectively; if $B$ is not a forcing set, then $\pt(G,B)=\infty$ or $\ptp(G,B)=\infty$.   %For $k\in\ZZ^+$, $\pt(G,k)=\min_{|S|=k}\pt(G,S)$ and $\ptp(G,k)=\min_{|S|=k}\ptp(G,S)$. %, these are called the \emph{$k$-standard propagation time of $G$} and \emph{$k$-PSD propagation time of $G$}, respectively. 
% The \emph{standard propagation time of $G$} (respectively, \emph{PSD propagation time of $G$}) is $\pt(G)=\pt(G,\Z(G))$ (respectively, $\ptp(G)=\ptp(G,\Zp(G))$). 
%When $B$ is a forcing set,  define the {\emph{round function}} by  $\rd_B(v)=k$ if and only if $v\in B^{(k)}$ for each $v \in V(G)$. 

 %  \lh{Forts moved from Section \ref{s:joinisfast}.}
    
    Let $G$ be a graph. A nonempty set  $F \subseteq V(G)$  is a standard fort of $G$ if there does not exist a vertex in $V(G) \setminus F$ that has exactly one neighbor in $F$.
%\end{definition}
%\begin{definition}
  Let $F \subseteq V(G)$ be nonempty, and $\{C_i\}_{i=1}^k$ be the components of $G[F]$. Then $F$ is a \emph{PSD fort} of $G$ if for each $v \in V(G) \setminus F$ and each $i \in \{1,2,\dots,k\}$, $|N_G(v) \cap V(C_i)| \neq 1$. 
%\end{definition}

\begin{theorem}
    Let $G$ be a graph and $B\subseteq V(G)$.
    \begin{enumerate}[$(1)$]
        \item {\rm\cite[Theorem 8]{BFH19}} The set $B$ is a zero forcing set of $G$ if and only if $B\cap F\ne \emptyset$ for every fort $F$ of $G$.
        \item {\rm {\cite[Theorem 2]{SMH19_PSDforts}}}  The set $B$ is a PSD forcing set of $G$ if and only if $B \cap F \neq \emptyset$ for each PSD fort $F$ of $G$.
    \end{enumerate}
\end{theorem}

%%%%%%%%%%%%%%%%%%%%%%%%%%%%%%%%%%%%%%%%%%%

%==========================
%
% Section
%
%==========================
\section{A   graph that does not have a full PSD propagation time interval}\label{s:PSDgap}
%\section{A family of graphs with gaps in the propagation time set and the expanded  propagation time set}\label{s:PSDgap}

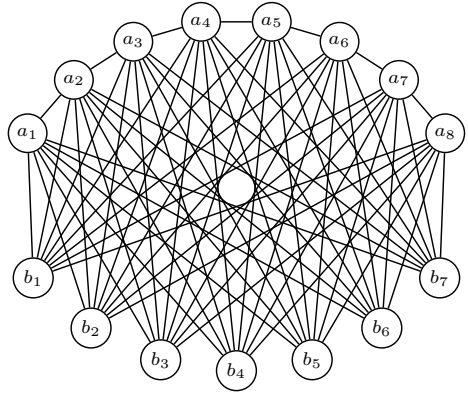
\begin{figure}[h!]\centering 
\begin{tikzpicture}[scale=.55,line width = .55\Lpt,vertex/.style={minimum size=20\Spt,draw,circle,node font=\scriptsize,inner sep=2pt},auto]
  % Left side: 7 vertices
  \foreach \i in {1,...,7}
    \node[vertex] (a\i) at  ($({-90+(\i-4)*180/11}:6.5)+(0,2)$) {$b_{\i}$};

  % Right side: 8 vertices
  \foreach \j in {1,...,8}
    \node[vertex] (b\j) at ($({90-(\j-4.5)*180/11}:6)+(0,-2)$) {$a_{\j}$};

  % Connect all edges between the two partitions
  \foreach \i in {1,...,7}
    \foreach \j in {1,...,8}
      \draw (a\i) -- (b\j);

  \draw (b1)--(b2)--(b3)--(b4)--(b5)--(b6)--(b7)--(b8);    
\end{tikzpicture}
%\scalebox{.5}{\includegraphics{c15}}
\caption{\label{f:gapex} The graph $\sgap 0$, which does not have a full PSD propagation interval nor a full  expanded PSD propagation time interval.}\end{figure}

 In this section we present with proof a family of graphs    whose PSD propagation time interval is not full (that includes Small's example), providing a negative resolution to Warnberg's conjecture in  \cite{PSDproptime}. We also show that the expanded PSD propagation time interval of each graph in the family is not full.  The smallest graph in the family    is shown in Figure \ref{f:gapex}.

\begin{definition}
  Let $A = \{a_1,\dots,a_{8+2k}\}$ and $B = \{b_1,\dots,b_{7+2k}\}$, and define a graph of order $15+4k$ by  $\sgap k= P_{8+2k}\vee (7+2k)K_1$ where $V(P_{8+2k})=A$ (in path order) and $V((7+2k)K_1)=B$. 
\end{definition}

%We show $\sgap k$ has gaps in the  PSD propagation time set and the expanded {\blu PSD} propagation time set.

%  and has \GraphSix"N]rEEB?^|~F{NxN|F|O"

\begin{theorem}\label{ExPSDproptimeGap}
   For $k\ge 0$, the graph $\sgap k$ of order $15+4k$   %shown in Figure \ref{f:gapex} 
   has $\Zp(\sgap k)=8+2k, \zpbar(\sgap k)=9+2k$, and $\EPTpLV(\sgap k)=\PTpLV(\sgap k)=\{1, 2, 4+k, 5+k, \dots,7+2k\}$. 
\end{theorem}

\bpf  We first note that any induced subgraph of $\sgap k$ containing at least one vertex each from $A$ and $B$ is connected. Let $X\subsetneq V(\sgap k)$ be a minimal PSD forcing set and define $W=V(\sgap k)\setminus X$. 

If $|W\cap A|\ge 2$ and $|W\cap B|\ge 2$, then $W$ is a PSD fort (because it induces a connected graph and any vertex in $X$ has two white neighbors in one of the sets $A$ or $B$) so $X$ is not a PSD forcing set.
Therefore $|X\cap A|\ge 7+2k$ or $|X\cap B|\ge 6+2k$.
We  analyze cases. 
\ben[$(i)$]
\item $|X\cap A|=8+2k$: The removal of $X$ results in  isolated vertices and $\ptp(\sgap k,X)=1$. Since $X$ is minimal, this implies $X\cap B=\emptyset$ and $|X|=8+2k$.

\item $|X\cap A|=7+2k$: 

If $X\cap B=\emptyset$, then $G[W]$ is connected and every vertex in $X$ has at least seven white neighbors in $B$, so no forcing could take place.

Thus $|X\cap B|\ge 1$.
Then  $(A\setminus\{a_i\})\cup \{b_j\}\subseteq X$ for some $i,j$, so in the first round $b_j$ forces $a_i$, and  $\ptp(\sgap k,X)=2$. Necessarily $|X\cap B|= 1$ and $|X|=8+2k$ since $X$ is minimal.  

\item $|X\cap B|=7+2k$: 

If $X\cap A=\emptyset$, then $G[W]$ is connected and every vertex in $X$ has at least eight white neighbors in $A$, so no forcing could take place.

Thus $|X\cap A|= 1$ and $|X|=8+2k$ because adding any one vertex of $A$ to $B$ produces a PSD forcing set. Observe that $\ptp(\sgap k,B\cup\{a_i\})=8+2k-i$ for $i=1,\dots,4+k$, yielding PSD propagation times $7+2k,6+2k,\dots,4+k$, respectively. Furthermore, $i\in\{5+k,\dots,8+2k\}$ yields the same set of PSD propagation times by symmetry. %Observe that $G[A]$ is a path and no vertex from $B$ can force a vertex in $A$ unless $X\cap B$ contains one of $a_1,a_2,a_7,a_8$.  Any one of the vertices added to seven vertices of $B$ results in a minimal PSD forcing set, with $\ptp(\sgap k,B\setminus\{b_i\}\cup\{a_1\})=$ In order for any Then $A\setminus\{a_i\}\du \{b_j\}\subseteq A$, in the first round $b_j$ forces $a_i$, and  $\ptp(\sgap k,X)=2$. Necessarily $|X\cap B|= 1$ since $X$ is minimal.

\item $|X\cap B|=6+2k$: 

If $|X\cap A|\le 1$, then $G[W]$ is connected, the blue vertex in $X\cap A$ (if there is one) has at least one white neighbor in $A$ and one white neighbor in $B$, and every vertex in $X\cap B$ has at least seven white neighbors in $A$, so  no forcing could take place.

Thus $|X\cap A|\ge 2$. Note first that the  vertices in $B$ are adjacent to every vertex in $A$, so as long  as there are two white vertices in $A$, no vertex in $B$ can force.  A vertex in $A$ cannot force a vertex in $A$ until the one white vertex in $B$ is forced.  In order for  vertex  $a_i$ to force a vertex in $B$, we must have $a_{i-1},a_i,a_{i+1}$ in $X$ (if all three vertices exist).  Thus we have two possibilities:

$|X\cap A|=2$ and $|X|=8+2k$.  In this case, $X\cap A=\{a_1,a_2\}$ or $X\cap A= \{a_7+2k,a_8+2k\}$ and $\pt(\sgap k,X)=7+2k$.

$|X\cap A|=3$ and $|X|=9+2k$.  In this case, $X\cap A=\{a_{i},a_{i+1},a_{i+2}\}$ for $i=2,\dots,5+2k$. %, $X\cap A=\{a_3,a_4,a_5\}$, $X\cap A=\{a_4,a_5,a_6\}$, or $X\cap A=\{a_5,a_6,a_7\}$.  
Let $X_B=X\cap B$.  Then $\ptp(\sgap k,X_B\cup\{a_{i},a_{i+1},a_{i+2}\})=7+2k-i$ for $i=2,\dots,3+k$, yielding PSD propagation times $5+2k,4+2k,\dots,4+k$, respectively. Furthermore, $i\in\{5+k,\dots,8+2k\}$ yields the same set of PSD propagation times by symmetry. 

%Then $A\setminus\{a_i\}\du \{b_j\}\subseteq A$, in the first round $b_j$ forces $a_i$, and  $\ptp(\sgap k,X)=2$. Necessarily $|X\cap B|= 1$ and $|X|=8$ since $X$ is minimal.  
\een

Thus  $\Zp(\sgap k)=8+2k$ and $\zpbar(\sgap k)=9+2k$, no minimal set has PSD propagation time  of $3,4,\dots,3+k$, $\EPTpLV(\sgap k)=\{1,2,4+k,5+k,\dots,7+2k\}$, and each of these values is realized by a minimum PSD forcing set, giving $\PTpLV(\sgap k)=\{1,2,4+k,5+k,\dots,7+2k\}$.\epf

%\begin{remark}Note that $\sgap k[A] \cong P_8$.  Constructing $\mathsf{S}_{15+2k}$ analogously to $\sgap k$ where $\mathsf{S}_{15+2k}[A]\cong P_{8+2k}$ (and $B$ as before), one derives a family of graphs with arbitrarily large gaps in the extended propagation time set, namely $\EPTpLV(\mathsf{S}_{15+2k})=\{1,2,4+k, \dots, 7+2k\}$.    \end{remark}

%%%%%%%%%%%%%%%%%%%%%%%%%%%%%%%%%%%%%%%%%%%

%==========================
%
% Section
%
%==========================

\section{Fixed propagation time greater than one}\label{s:bigfixedpt}

When studying the propagation time set or expanded propagation time set, it is natural to ask what  fixed propagation times are possible, in addition to which graphs have fixed propagation time.  In  Sections \ref{s:fastjoins} and \ref{s:joinisfast}  we discuss fixed propagation time equal to one. In this section we  present families of graphs that have larger standard fixed propagation time and show that no graph can have fixed PSD propagation time greater than one.

\subsection{Fixed standard propagation time greater than one}

We begin with examples of graphs that have fixed standard propagation time greater than one.

\begin{example}\label{ex:cycle}
   For the cycle $C_n$ with $n\ge 5$, a set is a minimal standard zero forcing set if and only if it consists of two adjacent vertices, and thus  $\PTul(C_n)=\PTbar(C_n)=\lc\frac{n-2}2\rc$.
\end{example}

\begin{example}\label{ex:star}
   For the star $K_{1,n-1}$ with $n\ge 3$, a set is a minimal standard zero forcing set if and only if it consists of $n-2$ leaves, and thus  $\PTul(K_{1,n-1})=\PTbar(K_{1,n-1})=2$.
\end{example}

The maximum size of a minimal standard forcing set for a graph $G$ is denoted by $\zbar(G)$. If $\Z(G) = \zbar(G)$, then every minimal standard forcing set of $G$ is also a minimum standard forcing set of $G$. Two minimum standard zero forcing sets $B_1$ and $B_2$ of a graph $G$ are \emph{isomorphic} if there is a graph automorphism $\psi$ of
$G$ such that $\psi(B_1)=B_2$.  Note that both of the previous examples are covered by the next observation.

\begin{observation}
    If $G$ is a graph with $\Z(G)=\zbar(G)$ and in which all minimum standard zero forcing sets are isomorphic, then $G$ has fixed standard propagation time.
\end{observation}

In the next examples we see examples of graphs where not all minimum standard zero forcing sets are isomorphic.  
\begin{example}\label{ex:order45}
   The smallest order $n$ for which there is a graph with fixed standard propagation time is greater than one is $n=3$, namely $K_{1,2}=P_3$. The connected graphs of orders $n=4$ and $n=5$ with fixed standard propagation time two are listed in Figure \ref{f:order4fixed}. %(with their standard propagation times).
    \begin{figure}[!ht]
 \begin{center}
  \scalebox{.5}{\begin{tikzpicture}[scale=.75,line width = 1.5\Lpt,vertex/.style={minimum size=30\Spt,draw,circle,thick},auto]
     \begin{scope}[shift={(2.5,0)}]
     \node[vertex] (1) at (0,0) {};
     \node[vertex] (2) at (90:2cm) {};
     \node[vertex] (3) at (210:2cm) {};
     \node[vertex] (4) at (330:2cm) {};
     \draw[thick] (2)--(1)--(3);
     \draw[thick] (1)--(4);
     \end{scope}
     
     \begin{scope}[shift={(6.75,.75)}]
     \node[vertex] (1) at (0:1.5cm) {};
     \node[vertex] (2) at (120:1.5cm) {};
     \node[vertex] (3) at (240:1.5cm) {};
     \node[vertex] (4) at (4,0) {};
     \draw[thick] (1)--(2)--(3)--(1)--(4);
     \end{scope}

     \begin{scope}[shift={(14,.5)}]
     \node[vertex] (1) at (45:2cm) {};
     \node[vertex] (2) at (135:2cm) {};
     \node[vertex] (3) at (225:2cm) {};
     \node[vertex] (4) at (-45:2cm) {};
     \draw[thick] (1)--(2)--(3)--(4)--(1)--(3);
     \end{scope} 

     \begin{scope}[shift={(0,-5)}]
     \node[vertex] (0) at (0,0) {};
     \node[vertex] (1) at (45:2cm) {};
     \node[vertex] (2) at (135:2cm) {};
     \node[vertex] (3) at (225:2cm) {};
     \node[vertex] (4) at (-45:2cm) {};
     \draw[thick] (1)--(0)--(2);
     \draw[thick] (3)--(0)--(4);
     \end{scope}
     
     \begin{scope}[shift={(5,-5)}]
     \node[vertex] (0) at (0,0) {};
     \node[vertex] (1) at (45:2cm) {};
     \node[vertex] (2) at (135:2cm) {};
     \node[vertex] (3) at (225:2cm) {};
     \node[vertex] (4) at (-45:2cm) {};
     \draw[thick] (4)--(1)--(0)--(2);
     \draw[thick] (3)--(0)--(4);
     \end{scope}

     \begin{scope}[shift={(10,-5)}]
     \node[vertex] (0) at (0,0) {};
     \node[vertex] (1) at (45:2cm) {};
     \node[vertex] (2) at (135:2cm) {};
     \node[vertex] (3) at (225:2cm) {};
     \node[vertex] (4) at (-45:2cm) {};
     \draw[thick] (4)--(1)--(0)--(2);
     \draw[thick] (3)--(0)--(4);
     \draw[thick] (4) --(3);
     \end{scope}

     \begin{scope}[shift={(15,-5)}]
     \node[vertex] (0) at (0,0) {};
     \node[vertex] (1) at (45:2cm) {};
     \node[vertex] (2) at (135:2cm) {};
     \node[vertex] (3) at (225:2cm) {};
     \node[vertex] (4) at (-45:2cm) {};
     \draw[thick] (4)--(1)--(0)--(2);
     \draw[thick] (3)--(0)--(4)--(3);
     \draw[thick,bend left] (4)edge node [bend left] {}(2);
     \end{scope}

     \begin{scope}[shift={(-2.5,-10)}]
     \node[vertex] (0) at (162:2cm) {};
     \node[vertex] (1) at (90:2cm) {};
     \node[vertex] (2) at (234:2cm) {};
     \node[vertex] (3) at (18:2cm) {};
     \node[vertex] (4) at (-54:2cm) {};
     \draw[thick] (1)--(0)--(3)--(1);
     \draw[thick] (3)--(4)--(2);
     %\draw[thick,bend left] (4)edge node [bend left] {}(2);
     \end{scope}

     \begin{scope}[shift={(2.5,-10)}]
     \node[vertex] (0) at (162:2cm) {};
     \node[vertex] (1) at (90:2cm) {};
     \node[vertex] (2) at (234:2cm) {};
     \node[vertex] (3) at (18:2cm) {};
     \node[vertex] (4) at (-54:2cm) {};
     \draw[thick] (0)--(2)--(4)--(3)--(1)--(0);
     %\draw[thick] (3)--(0)--(4);
     %\draw[thick,bend left] (4)edge node [bend left] {}(2);
     \end{scope}    

     \begin{scope}[shift={(7.5,-10)}]
     \node[vertex] (0) at (162:2cm) {};
     \node[vertex] (1) at (90:2cm) {};
     \node[vertex] (2) at (234:2cm) {};
     \node[vertex] (3) at (18:2cm) {};
     \node[vertex] (4) at (-54:2cm) {};
     \draw[thick] (2)--(4)--(3)--(1);
     \draw[thick] (3)--(2)--(1)--(4);
     \draw[thick] (0)--(1);
     %\draw[thick,bend left] (4)edge node [bend left] {}(2);
     \end{scope}

     \begin{scope}[shift={(12.5,-10)}]
     \node[vertex] (0) at (162:2cm) {};
     \node[vertex] (1) at (90:2cm) {};
     \node[vertex] (2) at (234:2cm) {};
     \node[vertex] (3) at (18:2cm) {};
     \node[vertex] (4) at (-54:2cm) {};
     \draw[thick] (2)--(4)--(3)--(1);
     \draw[thick] (3)--(0)--(4)--(1)--(2)--(3);
     %\draw[thick,bend left] (4)edge node [bend left] {}(2);
     \end{scope}

     \begin{scope}[shift={(17.5,-10)}]
             \node[vertex] (0) at (162:2cm) {};
     \node[vertex] (1) at (90:2cm) {};
     \node[vertex] (2) at (234:2cm) {};
     \node[vertex] (3) at (18:2cm) {};
     \node[vertex] (4) at (-54:2cm) {};
     \draw[thick] (0)--(2)--(4)--(3)--(1);
     \draw[thick] (3)--(0)--(4)--(1)--(2)--(3);
     %\draw[thick,bend left] (4)edge node [bend left] {}(2);
     \end{scope}
 \end{tikzpicture}
 }
 \caption{ The  connected graphs of order $n=4,5$  with $\PTul(G)=\PTbar(G)=2$ \label{f:order4fixed}}
 \vspace{-20pt}
 \end{center}
 \end{figure}
 
\end{example}

The star and all but two of the graphs in Figure \ref{f:order4fixed} are threshold graphs. %, which are a type of cograph.
A \emph{threshold graph} is %a graph with no induced $P_4$, $C_4$, or $2K_2$, or equivalently 
a graph that can be built from isolated vertices by a sequence of disjoint union and join operations, where at each stage one of the two graphs is a single isolated vertex. 
This can be diagrammed as a rooted binary tree where the leaves are isolated vertices and the vertices of degree two or more are marked as diamonds with $\du$ (for disjoint union) or $\vee$ (for join), with the restriction that at least one of the parts added in each operation is a single vertex, as illustrated in Figure \ref{fig:contree}; this is called a \emph{construction tree} of the graph. We show in Theorem \ref{p:threshfixed} that every threshold graph has fixed standard propagation time. Furthermore, every connected threshold graph except $K_n$ has fixed standard propagation time greater than one, as shown in Corollary \ref{ThreshPT1isKn}.

\begin{figure}[!ht]
 \begin{center}
  \begin{tikzpicture}[scale=0.7,construction graph settings]
     \node[vertex] (1) at (-90:4cm) {$1$};
     \node[vertex] (2) at (40:4cm) {$2$};
     \node[vertex] (3) at (140:4cm) {$3$};
     \node[vertex] (4) at (-.75,0) {$4$};
     \node[vertex] (5) at (.75,0) {$5$};
     \draw (4)--(5)--(3)--(4);
     \draw (5)--(2)--(4);
     \draw (5)--(1)--(4);

     \begin{scope}[shift={(4,-4)}]
     \node[vertex] (1) at (0,0) {$1$};
     \node[vertex] (2) at (1.5,0) {$2$};
     \node[vertex] (3) at (3,1.5) {$3$};
     \node[vertex] (4) at (4.5,3) {$4$};
     \node[vertex] (5) at (6,4.5) {$5$};

     \node[operation] (6) at (1.5,1.5) {$\du$};
     \node[operation] (7) at (3,3) {$\du$};
     \node[operation] (8) at (4.5,4.5) {$\vee$};
     \node[operation] (9) at (6,6) {$\vee$};

     \begin{scope}[construction edges]
     \draw (1)--(6)--(7)--(8)--(9)--(5);
     \draw (2)--(6);
     \draw (3)--(7);
     \draw (4)--(8);
     \draw (5)--(9);
     \end{scope}
         
     \end{scope}
 \end{tikzpicture}
 \caption{  A threshold graph and its construction tree\label{fig:contree}}
\vspace{-10pt}
 \end{center}
 \end{figure}

We begin with elementary results on the standard zero forcing number of graphs with a universal vertex, a class that includes connected threshold graphs.  These results are known and have been used elsewhere (e.g., \cite{minimalzfsets}), but we do not know of  clear sources for some.
\begin{proposition}\label{prop:universal}
    Let $G$ be a graph and let  $H=G\vee K_1$. 
    \begin{enumerate}[$(1)$]
        \item\label{l1}{\rm \cite[Proposition 9.16]{HLS22book}} If $G$ has no isolated vertices, then $\Z(H) = \Z(G)+1$.
        \item\label{l3} If $G$ has at least two isolated vertices, then $\Z(H) = \Z(G)-1$.
        \item\label{l2} If $G$  has exactly one isolated vertex, then $\Z(H) = \Z(G)$.    \end{enumerate}
    \end{proposition}
    \begin{proof}
       %  Item \eqref{l1} appears  in \cite[Proposition 9.16]{HLS22book}.
         Let  $V(K_1)=\{u\}$.
         For Item \eqref{l3}, let $I=\{v_1,v_2,\dots,v_k\}$ denote the isolated vertices of $G$, and assume that $k\geq 2$.   Choose a minimum standard zero forcing set $B$ for $G$, which necessarily contains $I$ and let $\mathcal{F}$ be a chronological list of forces for $B$ in $G$. Let $B^{\prime} = B\setminus \{v_1\}$.  Then  $v_2 \to u$ is a valid first force for $B^{\prime}$ in $G$, after which all of the forces of $\mathcal{F}$ are valid in $H$, and only $v_1$ remains white.  Adding the final force $u \to v_1$ shows that $B^{\prime}$ is a standard zero forcing set of $H$ and thus $\Z(H) \leq \Z(G)-1$, and it is well known that $\Z(H)\ge \Z(G)-1$ since $G=H-u$ \cite[Proposition 9.13]{HLS22book} %\cite{EHHLR12}.
         
         To show Item \eqref{l2}, first  let $B$ be a minimum standard zero forcing set of $G$, $\mathcal{F}$ be a chronological list of forces for $B$ in $G$, and $v$ be the isolated vertex of $G$.  Then $v\in B$, and since $u$ is the only neighbor of $v$ in $H$, the force $v \to u$ is a valid first force in $H$ starting with the vertices in $B$ blue.  After this force, every force in $\mathcal{F}$ is valid in $H$, so $B$ is a standard zero forcing set of $H$.  Thus $\Z(H)\leq \Z(G)$. 

         Now suppose $B'$ is a minimum standard zero forcing set of $H$ with corresponding chronological list of forces $ \mathcal{F}$.  If $v \in B'$, then  $u \not\in B'$ since $B'$ is minimal and  the first force must be $v \to u$; the remaining forces of $ \mathcal{F}$  can take place entirely in $G$  (replacing a last force $u\to x$ by $y\to x$ for $y\in N_G(x)$ if necessary).  Thus %\st{deleting $v \to u$ from $\mathcal{F}$ yields a chronological list of forces for $B$ in $H$ that results in vertices of $H$ blue, so} 
         $\Z(G) \leq \Z(H)$.  If $v \not\in B'$ and $u \not\in B'$, then the first force in $\mathcal{F}$ must be $x \to u$, where $x\neq v$.  Since $u$ is the unique white neighbor of $x$ in $H$, we have $N_G(x) \subseteq B'$.  Thus we may obtain a new minimum standard zero forcing set $ (B' \setminus \{y\}) \cup \{u\}$, where $y\in N_G(x)$.  Therefore if $v \not\in B'$ we may assume $u \in B'$. Under this hypothesis, every force in $\mathcal{F}$ is valid in $G$ except $u \to v$, and this must be the last force in $\mathcal{F}$. Thus $ (B' \setminus \{u\} )\cup \{v\}$  is a standard  zero forcing set of $G$ whose chronological list of forces is $\mathcal{F}$ except the last force, so % and let $\mathcal{F}^{\prime}$ be the chronological list of forces obtained by deleting $u\to v$ from $\mathcal{F}$. 
         %Then $\mathcal{F}^{\prime}$ is a chronological list of forces for $B^{\prime}$ in $H$ that colors all vertices of $H$ blue, and therefore 
         $\Z(G) \leq \Z(H)$.
    \end{proof}

%\lh{upgraded both \Cref{p:minimalfixed} and \Cref{p:threshfixed} from props to thms}

\begin{theorem}\label{p:minimalfixed}
    If $G$ is a  threshold graph, then $\Z(G)=\zbar(G)$.
\end{theorem}
 \bpf Let $G$ be a threshold graph with $|G|=n$.  We proceed by induction on $n$.  The statement is true if $n\leq 2$, so suppose the statement holds for $k\leq n$.  Let $H$ be a threshold graph with $|H|=n+1$. The statement is trivial if $H$ is the form $G\du K_1$, so consider $H=G\vee K_1$, where $G$ is a threshold graph. Let $u$ denote the universal vertex joined to $G$.
Let $B$ be a  minimal standard zero forcing set of $H$ of maximum cardinality. 

If $G$ has two or more isolated vertices, then $\Z(H) = \Z(G) -1$  by Proposition \ref{prop:universal}\eqref{l3} %, \Cref{l3} \Cref{prop:universal}, \Cref{l3} } 
and $B$ contains all but one of the isolated vertices of $G$. Furthermore, $u \not\in B$ since any of the isolated vertices of $G$ in $B$ can force $u$  of $H$ in the first round of propagation.   Let $x$ be the isolated vertex in $G$ such that $x\not\in B$.  Then $B^{\prime} = B \cup \{x\}$ is a standard zero forcing set of $G$.  Suppose $B^{\prime}$ is not minimal for $G$ and $y$ is a vertex such that $B^{\prime} \setminus \{y\}$ is a standard zero forcing set of $G$.  Then $y$ cannot be any of the isolated vertices of $G$, and consider $B\setminus \{y\}$.  Let $w \neq x$ be an isolated vertex of $G$.  Then $w \rightarrow u$ is a valid first force for $B\setminus \{y\}$ in $H$, and each of the forces occurring in a chronological list of forces for $B^{\prime} \setminus \{y\}$ in $G$ is valid for $B \setminus \{y\}$ in $H$.  Thus $B \setminus \{y\}$ is a standard zero forcing set of $H$, contradicting the minimality of $H$. Therefore, $\zbar(H) = |B| = |B^{\prime}| -1 = \Z(G) -1 = \Z(H)$.

{So we assume $G$ has either one or no isolated vertices. 
Suppose first that $G$ has one isolated vertex $x$ and $x\in B$. 
Then $u\not\in B$ and the only force in the first round is $x\to u$. 
Let $B'=(B\setminus\{x\})\cup \{u\}$ and note that $B'$ is a standard forcing set of the same cardinality because $B'$ uses the forces of $B$ except that the first round force $x\to u$ is replaced by a  last round force $u\to x$.  
So henceforth we assume no isolated vertex of $G$ is in $B$.}

If $u \not\in B$, then the first force in $H$ must be of the form $z \rightarrow u$, where $z\in V(G)$, $N_G[z] \subseteq B$, and there exists  $y \in N_G(z) \cap B$.  
Then $B^{\prime}=(B \setminus \{y\}) \cup \{u\}$ is a minimal standard zero forcing set of $H$.  Thus we assume for the remainder of the proof that $u \in B$.

Now suppose $G$ has no isolated vertices (whether in $B$ or not). Then $\Z(H) = \Z(G) +1$  by Proposition \ref{prop:universal}\eqref{l1}, and $B^{\prime}=B\setminus \{u\}$ is a standard zero forcing set of $G$ since $u\rightarrow w$ can only occur as the last force in $H$, which can be replaced $y \rightarrow w$ for any $y \in N_G(w)$. 
We claim that $B^{\prime}$ is a minimal standard zero forcing set of $G$. Suppose to the contrary that 
 $B'\setminus \{y\}$ is a standard zero forcing set for $G$, where $y \in B \cap V(G)$.  Then $B\setminus \{y\}$ is a standard  zero forcing set for $H$ since each of the forces for $B^{\prime} \setminus \{y\}$ in $G$ is valid for $B \setminus \{y\}$ in $H$, contradicting the minimality of $H$.  This gives $\zbar(H) = |B| = \Z(G) + 1 = \Z(H)$.

Finally, suppose $G$ has exactly one isolated vertex $x$.  Then $\Z(H)=\Z(G)$  by Proposition \ref{prop:universal}\eqref{l2}, and since  $u\in B$ it must be the case that $x \not\in B$ since $B$ is minimal and $u \rightarrow x$ is the final force performed.  Consider $B^{\prime} = (B \setminus \{u\}) \cup \{x\}$, which is a standard zero forcing set of $H$, as can be seen by adding $x\rightarrow u$ as the first force, proceeding with the forces for $B$ in $H$ until the final force $u\to x$, which is removed. 
As before, if $B^{\prime} \setminus \{y\}$ is a standard zero forcing set of $G$, then $B \setminus \{y\}$ is a standard zero forcing set of $H$, a contradiction {(note that $y\ne x$ since $x$ is an isolated vertex in $G$)}.  This gives $\zbar(H) = |B| = |B^{\prime}| = \Z(G) = \Z(H)$.
\epf

 The previous theorem combined with the next result establishes threshold graphs have fixed propagation time. The  \emph{(standard) throttling number}  of $B$ in $G$ is $\thz(G,B)=|B|+\pt(G,B)$ and the    \emph{(standard) throttling number} of  a graph $G$  is  $\thz(G)=\min_{B\subseteq V(G)} \thz(G,B) $.   Observe that $\thz(G,B)\le |V(G)|$ for any standard zero forcing $B$, because at least one vertex must be forced in each round.

\begin{theorem}{\rm \label{t:thz-hi-ext}\cite[Theorem 11.22]{HLS22book}}  A graph $G$ has $\thz(G)=|V(G)|$
if and only if $G$ is a threshold graph.
\end{theorem}

\begin{corollary}
\label{p:threshfixed}
    Every threshold graph $G$ has fixed standard propagation time equal to $|V(G)|-\Z(G)$.
\end{corollary}

\begin{proof}
    If $B$ is a minimal standard zero forcing set of a threshold graph $G$, then $|B|=\Z(G)=\zbar(G)$ by Theorem \ref{p:minimalfixed}.  By Theorem \ref{t:thz-hi-ext}, $\thz(G)=|V(G)|$, so $\thz(G,B)=|V(G)|$ for every minimum standard zero forcing set. Thus $\pt(G,B)=|V(G)|-|B|=|V(G)|-\Z(G)$.
\end{proof}

 The next result is immediate from the previous corollary.
\begin{corollary}
\label{ThreshPT1isKn}
If $G$ is a connected threshold graph of order $n\ge 2$ with $\PTbar(G)=1$, then $G\cong K_n$.
\end{corollary}

A \emph{cograph} is %a graph with no induced $P_4$, or equivalently 
a graph that can be built from isolated vertices by a sequence of disjoint union and join operations (the restriction that one of the two graphs in each operation is a single vertex is removed).  
There are cographs that do not have fixed standard propagation time and for which $\Z(G)<\zbar(G)$, as seen in the next two examples. 

\begin{figure}[!ht]
      \begin{center}
      \begin{tikzpicture}[scale=0.8,construction graph settings]
      \begin{scope}[scale=2]

          \node[vertex] (1) at (0,0) {$1$};
          \node[vertex] (2) at (2,0) {$2$};
          \node[vertex] (3) at (2,-2) {$3$};
          \node[vertex] (4) at (0,-2) {$4$};
          \node[vertex] (5) at (1,-1) {$5$};
          \draw (1)--(2)--(3)--(4)--(1)--(5)--(2);
          \draw (4)--(5)--(3);
    \end{scope}
          \begin{scope}[shift={(6,-4.25)}]
              \node[vertex] (a) at (0,0) {$1$};
              \node[vertex] (b) at (1.25,0) {$2$};
              \node[vertex] (c) at (2.5,0) {$3$};
              \node[vertex] (d) at (3.75,0) {$4$};

              \node[operation] (e) at (1.25,1.5) {$\du$};
              \node[operation] (f) at (2.5,1.5) {$\du$};

              \node[operation] (g) at (2.5,3) {$\vee$};
              \node[vertex] (h) at (3.75,3) {$5$};
              \node[operation] (i) at (3.75,4.5) {$\vee$};

              \begin{scope}[construction edges]
              \draw (a)--(e) --(g)--(i)--(h);
              \draw (c) --(f)--(g);
              \draw (b)--(e);
              \draw (d)--(f);
              \end{scope}
          \end{scope}
      \end{tikzpicture}
 \caption{ The cograph $W_5$ and a construction tree\label{fig:contree-W5}}

 \end{center}
 \end{figure}
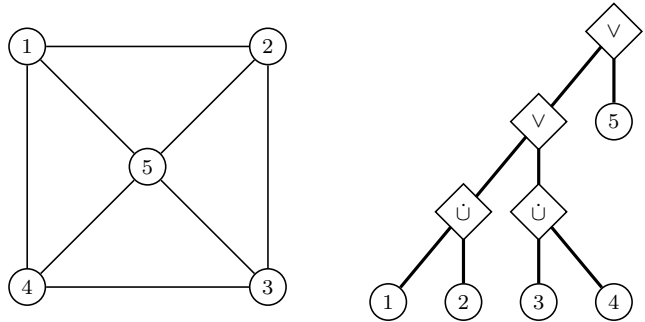

\begin{example}
     The  wheel graph $W_5$ is shown in Figure \ref{fig:contree-W5} together with a construction tree.  Note that $W_5$ does not have fixed standard propagation time because (with the labeling shown), $\pt(W_5,\{1,2,5\})=1$ and  $\pt(W_5,\{1,2,3\})=2$.
     \end{example} 

\begin{figure}[!ht]
    \centering
    \begin{tikzpicture}[scale=.5,construction graph settings]

    \node[vertex] (1) at (0,0) {$1$};
    \node[vertex] (2) at (4,0) {$2$};
    \node[vertex] (3) at (0,6) {$3$};
    \node[vertex] (4) at (4,6) {$4$};
    \node[vertex] (5) at (0,3) {$5$};
    \node[vertex] (6) at (4,3) {$6$};
    \node[vertex] (7) at (2,6) {$7$};
    \draw (1)--(2)--(6);
    \draw (5)--(1)--(6)--(4)--(5)--(3)--(6);
    \draw (2)--(5);
    \draw (7)--(1);
    \draw (7)--(2);
    \draw (7)--(3);
    \draw (7)--(4);
    \draw (7)--(5);
    \draw (7)--(6);

    \begin{scope}[shift={(6,0)}]
       \node[vertex] (1) at (0,0) {$1$};
       \node[vertex] (2) at (1.5,0) {$2$};
       \node[vertex] (3) at (3,0) {$3$};
       \node[vertex] (4) at (4.5,0) {$4$};
       \node[vertex] (5) at (6,0) {$5$};
       \node[vertex] (6) at (7.5,0) {$6$};
       \node[vertex] (7) at (9,0) {$7$};

       \node[operation] (a) at (.75,1.5) {$\vee$};
       \node[operation] (b) at (3.75,1.5) {$\du$};
       \node[operation] (c) at (6.75,1.5) {$\du$};
       \node[operation] (d) at (9,6) {$\vee$};
       \node[operation] (e) at (2.25,3) {$\du$};
       \node[operation] (f) at (5.25,6) {$\vee$};
       
       \begin{scope}[construction edges]
       \draw (1)--(a)--(2);
       \draw (3)--(b)--(4);
       \draw (5)--(c)--(6);
       \draw (a)--(e)--(b);
       \draw (e)--(f)--(c);
       \draw (f) -- (d) -- (7);
       \end{scope}
    \end{scope}
        
    \end{tikzpicture}
    \caption{A cograph $G$ with $\Z(G)\ne \zbar(G)$ and  a construction tree}
    \label{fig:contree-Z-Zbar}
\end{figure}
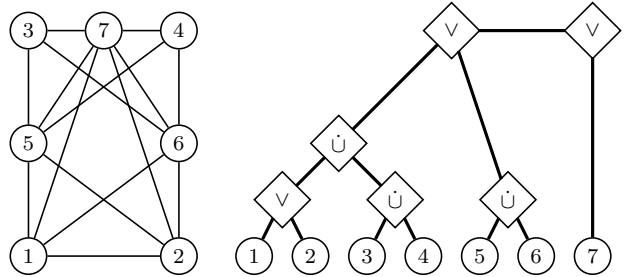

\begin{example}
     The  cograph $G$  shown in Figure \ref{fig:contree-Z-Zbar} (together with a construction tree) has $Z(G)=4<5\le\zbar(G)$ because $\{1,3,5,7\}$  and $\{1, 2, 3, 4, 5\}$ are minimal standard  zero forcing sets.  Any standard zero forcing set must contain one of the vertices 1 and 2, one of the vertices 3 and 4, and one of the vertices 5 and 6 because each of these pairs are twins \cite[Proposition 9.15]{HLS22book}. Furthermore, $\{1,3,5\}$  is not a standard zero forcing set. (In fact, $\zbar(G)=5$ as verified in \cite{sage}.)
     \end{example}

Finally we note that even if a graph $G$ satisfies $\zbar(G) \neq \Z(G)$, it may still be the case $G$ has fixed standard propagation time.   Two minimal standard zero forcing sets $B_1$ and $B_2$  of a graph $G$ are \emph{isomorphic} if there is a graph automorphism $\psi$ of
$G$ such that $\psi(B_1)=B_2$ (necessarily $|B_1|=|B_2|$).
 
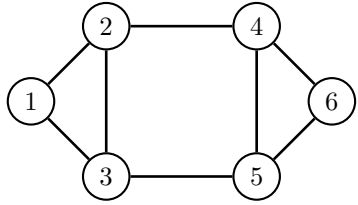
\begin{figure}[!ht]
    \centering
    \begin{tikzpicture}[scale=1,line width = 1\Lpt,vertex/.style={minimum size=30\Spt,draw,circle,thick},auto]

    \node[vertex] (1) at (0,0) {$1$};
    \node[vertex] (2) at (1,1) {$2$};
    \node[vertex] (3) at (1,-1) {$3$};
    \node[vertex] (4) at (3,1) {$4$};
    \node[vertex] (5) at (3,-1) {$5$};
    \node[vertex] (6) at (4,0) {$6$};
    \draw (1)--(2)--(3)--(1);
    \draw (2)--(4)--(6)--(5)--(4);
    \draw(3)--(5);
        
    \end{tikzpicture}
    \caption{A graph $G$ with $\Z(G)\neq \zbar(G)$ and fixed standard propagation time}
    \label{f:znezbar}
\end{figure}

\begin{example}\label{e:znezbar}
Let $G$ be the graph in Figure \ref{f:znezbar}.  Then, every minimum standard forcing set is isomorphic to $B_2=\{ 1,2\}$ because such a set must contain a vertex of degree two and a neighbor. Every minimal standard forcing set of size three is  isomorphic to $B_3=\{2,3,4\}$ because it can't contain a vertex of degree two and must contain a vertex of degree three and two of its neighbors.  Moreover, no set of four or more vertices is a minimal forcing set.  Direct computation reveals that $\pt(G,B_2)=3=\pt(G,B_3)$. %\st{Then, up to graph isomorphism, the only minimum standard forcing set is $\{ 1,2\}$. There is (up to graph isomorphism) one additional minimal standard forcing set $\{2,3,4\}$.  Direct computation reveals that $\PTbar(G) = \pt(G)$=3.}
    
\end{example}

%\newpage
\subsection{No graph has fixed PSD propagation time greater than one}
We  apply  known results to show that the only possible fixed PSD propagation time is one.

\begin{lemma}\label{L:PSD_G-Sconn}{\rm \cite[Lemma 4]{DSS24}}
If $G$ is a connected graph, then there exists a minimum positive semidefinite forcing set $B \subseteq V(G)$ such that $G - B$ is connected.
\end{lemma}

\begin{lemma}\label{L:PSD_ptp_reduce}{\rm \cite[Lemma 3.1]{MRC-PSDproptime}}
Let $G$ be a graph with PSD forcing set $B$ such that $G-B$ is connected.
Let $B'$ be the endpoints of the PSD forcing trees after the first time step. Then $B'$ is another PSD forcing set of $G$ with $|B|=|B'|$. Furthermore, if $\ptp(G,B)\geq 2$, then $\ptp(G,B')=\ptp(G,B)-1$.
\end{lemma}

\begin{theorem}
Let $G$ be a connected graph.  If $G$ has fixed PSD propagation time, then $\ptp(G)=1$.
\end{theorem}
 \bpf 
Suppose $G$ has a minimum PSD forcing set $B$ such that $B$ is a minimum forcing set with $\ptp(G,B)\ge 2$.  By Lemma \ref{L:PSD_G-Sconn}, there exists a minimum PSD forcing set $\hat B$ such that $G - \hat B$ is connected.  If $\ptp(G,\hat B)\ne \ptp(G,B)$, then $G$ does not have fixed PSD propagation time.  So assume $\ptp(G;\hat B)= \ptp(G,B)\ge 2$.  Then by Lemma \ref{L:PSD_ptp_reduce} there is another minimum PSD forcing set $B'$ such that $\ptp(G,B')<\ptp(G;\hat B)$.  This implies $G$ does not have fixed PSD propagation time.
\epf

%%%%%%%%%%%%%%%%%%%%%%%%%%%%%%%%%%%%%%%%%%%

%==========================
%
% Section
%
%==========================

\section{Fast joins have fixed propagation time equal to one}\label{s:fastjoins}
In this section we show that fast joins have fixed propagation time equal to one.  We also show that fast joins have very large forcing numbers.   We begin with  some preliminary material.

\begin{observation}
    A graph $G\not\cong K_n$ is a PSD fast join if and only if $\Gc\cong \du_{i=1}^t K_{n_i,m_i}$ where $n_i,m_i\ge 1$ for $i=1,\dots,t$  and $t\ge 2$. A graph $G\not\cong K_n$ is a standard fast join if and only if $\Gc\cong \du_{i=1}^t G_i$ where $t\ge 2$ and for each $i=1,\dots,t$,  $G_i\cong K_{n_i}$ %with $n_i\ge 1$
    or $G\cong K_{n_i,m_i}$ with $n_i,m_i\ge 2$.
\end{observation}

%\bs{Consider ways to rewrite this in smaller pieces: signatures? forts? lemmas?}

%\bs{Introduce forts; results on duplication as motivation.}
  
%\begin{definition}%[\cite{stdfort}\label{defstdfort}]

\begin{observation}\label{twinforts}
    Let $G$ be a graph.   
    \begin{enumerate}
        \item If $v_1, v_2 \in V(G)$ satisfy $N_G[v_1] = N_G[v_2]$ (implying $v_1v_2 \in E(G)$), then $\{v_1,v_2\}$ is a PSD fort and therefore a standard fort.
        \item If $v_1, v_2 \in V(G)$ satisfy $N_G(v_1) = N_G(v_2)$ (implying $v_1v_2 \not\in E(G)$), then $\{v_1,v_2\}$ is a standard fort.
    \end{enumerate}
\end{observation}

\begin{lemma}\label{componentforts}
    Let $G = \bigvee_{i=1}^k G_i$, let $B$ be a minimum (standard or PSD) forcing set of $G$, and let $W= V(G) \setminus B$. Then  
         $W\cap V(G_i) \neq \emptyset$ for at most two values of $i \in \{ 1,2, \dots, k \}$.
\end{lemma}
\begin{proof}
    %To show \eqref{f1}, s
    Suppose to the contrary  and without loss of generality that $|W \cap V(G_i)| \geq 1$ for $i \in \{1,2,3\}$, and choose $v_i \in W \cap V(G_i)$ for each $i$.  Then $G[\{v_1,v_2,v_3\}] \cong K_3$, and any $b \in B$ is adjacent to at least two vertices in $\{v_1,v_2,v_3\}$.  Thus no forces (standard or PSD) are possible, so $B$ is not a (standard or PSD) forcing set.  %The proof for \eqref{f2} is similar. \lh{NOT}
\end{proof}

% \begin{definition}
% Let $G = \bigvee_{i=1}^k G_i$ and $B$ a minimum (standard, PSD) zero forcing set of $G$.  The join signature of $B$, $s(B)$, is defined by \[s(B) = (a_1,a_2,\dots, a_k) \] where $a_i = \vert B^C \cap V(G_i) \vert$.
% \end{definition}

% \begin{definition}
% Given a join signature $s(B)=(a_1,a_2,\dots, a_k)$, a sequence $\sigma(B)$ is a subsignature of $s(B)$ if $\sigma(B) = (a_{j_1}, a_{j_2}, \dots , a_{j_r})$ where each $j_s \in \{ 1, 2, \dots, k\}$.
% \end{definition}

% \begin{lemma}
% \label{fort_signatures}
% If $B$ is a (standard zero, psd) forcing set for $G = \bigvee_{i=1}^k G_i$, then $s(B)$ has no subsignatures of the form 
% \begin{enumerate}
%     \item $(a_1, a_2)$ such that $a_1,a_2 \geq 2$
%     \item $(a_1,a_2,a_3)$ such that $a_1, a_2, a_3 \geq 1$ 
% \end{enumerate}
% \end{lemma}

\begin{theorem}
\label{t:ConnPSDClassHasFPT1}\label{t:ConnStdClassHasFPT1} 
If $G$ is a PSD fast join, then $\PTpbar(G)=1$  and $\Zp(G)\ge |V(G)|-2$.  If $G$ is a standard fast join, then $\PTbar(G)=1$  and $\Z(G)\ge |V(G)|-2$.
\end{theorem}

\begin{proof}
Assume $G$ is a PSD fast join or a standard fast join of order $n$.  
If $G$ is a complete graph of order at least two, then $\Zp(G)=\Z(G)=n-1$, and so $\PTpbar(G)=\PTbar(G)=1$.  
So assume  $G\not\cong K_n$, which implies $\Zp(G),\Z(G)\le n-2$. %$\Gc$ has at least two components, 
Let $B$ be a  minimal (standard or PSD) forcing set of $G$ and  let $W=V(G)\setminus B$, so that $|B|\le n-2$ and $|W|\ge 2$. 

  In light of Lemma \ref{componentforts}, $W$ can contain
vertices from at most two components of $\overline{G}$. By Observation \ref{twinforts},  we may assume $G$ has no  white adjacent twins.  Thus, $W$ can contain at most one vertex from each partite set of a component of $\overline{G}$ that is (complete) bipartite. 
Observation \ref{twinforts} also allows us to assume no  white independent twins if $G$ is a standard fast join; thus $W$ can contain at most one vertex from a component of $\overline{G}$ that is a complete graph ($K_2=K_{1,1}$ is treated as a complete graph for a standard fast join).
  We are thus left with the following possibilities for $W$:
\begin{enumerate}
%\item $W=\emptyset$, in which case $\pt_+(G;S)=0$;
%\item $|W|=1$;
\item\label{case3:ConnPSDClassHasFPT1} $W$ consists of one vertex from each of two components of $\overline{G}$;
\item\label{case4:ConnPSDClassHasFPT1} $W$ consists of one vertex from each partite set of a single component of $\overline{G}$; or 
%\item\label{case5:ConnPSDClassHasFPT1} $W$ consists of one vertex from each partite set of one component of $\overline{G}$ and one vertex from another component of $\overline{G}$; or
%\item\label{case6:ConnPSDClassHasFPT1} $W$ consists of one vertex from each partite set of two components of $\overline{G}$.
\item\label{case5:ConnPSDClassHasFPT1} $W$ consists of one vertex from each partite set of one component of $\overline{G}$ and  at least one vertex from another component of $\overline{G}$. 
\end{enumerate}
%Note that the white vertices induce a connected graph in $G$, so only the standard zero forcing rule applies  initially, in all cases.

First we show that the set $B$ in case \eqref{case5:ConnPSDClassHasFPT1} %and \eqref{case6:ConnPSDClassHasFPT1} 
is neither a PSD nor a standard forcing set. Note that in this case, the white vertices induce a connected graph in $G$, so only the standard zero forcing rule applies  initially. 

\begin{itemize}
\item Any blue vertex in a component of $\overline{G}$ with two white vertices is adjacent in $G$ to the white vertex in its partite set and a white vertex in the other component of $\overline{G}$ with white vertices,
\item Any blue vertex in a component of $\overline{G}$ with one white vertex is adjacent in $G$ to both of the white vertices in the component with two white vertices, and
\item Any other blue vertex is adjacent to all of the white vertices.
\end{itemize}
Because $B$ cannot perform any forces, it is not a (standard or PSD) forcing set.  Since case \eqref{case5:ConnPSDClassHasFPT1} 
cannot occur  for PSD or standard forcing and $B$ is any minimal forcing set,  the (standard or PSD) forcing number is at least $ n-2$.

For cases \eqref{case3:ConnPSDClassHasFPT1} and \eqref{case4:ConnPSDClassHasFPT1}, we consider PSD forcing and standard forcing separately, beginning with PSD forcing.   
 For case \eqref{case4:ConnPSDClassHasFPT1}, both white vertices can be PSD forced in the first round by any vertex in any of the other components, because the two white vertices are in separate components of $G-B$. Now consider case \eqref{case3:ConnPSDClassHasFPT1}. For $i=1,2$, let $A_i$ (containing white vertex $w_i$) and $B_i$ (all vertices blue)  denote the partite sets of the two components  of $\overline{G}$ containing white vertices. Then in the first round, $w_1$ can be forced  by a vertex in $B_2$ and  $w_2$ can be forced  by a vertex in $B_1$.    
 In both cases, $\ptp(G,B)=1$. 
   Thus, for all PSD forcing sets $B\subsetneq V(G)$ of $G$, we have $\ptp(G,B)=1$ and therefore $\PTpbar(G)=1$.  
   
   Finally we consider standard forcing.  For case \eqref{case4:ConnPSDClassHasFPT1}, both white vertices can be standard forced in the first round, each by a blue vertex in the same partite set of the component in $\Gc$, because each partite set contains at least two vertices.
   For case \eqref{case3:ConnPSDClassHasFPT1}, denote  the vertices of the components of $\Gc$ that contain white vertices by $U_1$ and $U_2$ and the white vertices by $w_i$ with $w_i\in U_i$ for $i=1,2$.  Then $w_1$  can be forced in round one by a blue vertex in $U_2$  that is not adjacent in $G$  to $w_2$, and similarly for $w_2$.   % there are several possibilities: let the two components that contain white vertices have partite sets $A_i$ (containing white vertex $w_i$) and $B_i$ (all vertices blue) for $i=1,2$. Then in the first round, $w_1$ can be forced  by a vertex in $B_2$ and  $w_2$ can be forced  by a vertex in $B_1$.    
  In both cases, $\pt(G,B)=1$. 
   Thus, for all standard zero  forcing sets $B\subsetneq V(G)$ of $G$, we have $\pt(G,B)=1$ and  therefore $\PTbar(G)=1$.   
    %\end{proof}
%
%
%
%\begin{corollary}If $G$ is a PSD fast join, then $\Zp(G)\geq |V(G)|-2$.  If $G$ is a standard fast join, then $\Z(G)\geq |V(G)|-2$. \end{corollary}
%
%\begin{proof}
\end{proof}

%%%%%%%%%%%%%%%%%%%%%%%%%%%%%%%%%%%%%%%%%%%
%==========================
%
% Section
%
%==========================.

\section{Joins with fixed propagation time equal to one are fast joins}
\label{s:joinisfast}

 In this section we resolve Conjectures \ref{c:psd} and \ref{c:std} under the additional hypotheses that the graph is a join of two graphs, or that the forcing number is sufficiently large. We begin with a few necessary tools. 

\begin{lemma} \label{inc-neigh}
Let $G$ be a %connected 
graph and 
 let $u$ and $v$ be adjacent vertices of $G$ such that $N_G[u]\subsetneq N_G[v]$. Then $\pt(G, V(G)\setminus \{u,v\})=\ptp(G, V(G)\setminus \{u,v\})= 2$. 
\end{lemma}

\begin{proof} %(a) 
Start with exactly the vertices in $B=V(G)\setminus\{u,v\}$ blue. Then no vertex can force $u$ in the first round but there is a vertex $x$ adjacent to $v$ and not to $u$, so $x\to v$ in the first round and $v\to u$ in the second round.  Thus $\pt(G,B)=\ptp(G,B)=2$.
\end{proof}

The next two corollaries are immediate from the previous lemma.
\begin{corollary}
 If $G$ is a connected graph of order $n\ge 3$ with $\PTpbar(G)=1$ or $\PTbar(G)=1$, then $G$ has no leaves.
\end{corollary}
%\begin{proof} In a connected graph of order at least three, a leaf is properly cornered by its neighbor.\end{proof}

\begin{corollary}
\label{NoDomVtx}
If $G\not\cong K_n$ and $\PTpbar(G)=1$ or $\PTbar(G)=1$, then $G$ cannot have a universal vertex.
\end{corollary}

%\begin{proof}A universal vertex corners all its neighbors.  If $G\neq K_n$ then some neighbor is properly cornered.\end{proof}

A set $B\subsetneq V(G)$ is a \emph{slow} standard (respectively, PSD) forcing set if $\pt(G,B)\ge 2$ (respectively, $\ptp(G,B)\ge 2$).

\begin{lemma}
\label{Join2ForbiddenInducedSG}
If $G$ is a join of two  graphs and at least one  of the two conditions
$ \PTpbar(G)=1$ or $\PTbar(G)=1$ holds, then none of the graphs $P_4$, $P_3 \du K_1$,
or $K_2 \du 2K_1$ is an induced subgraph of $G$.
\end{lemma}

\begin{figure}[!ht]
\centering
\scalebox{.6}{\begin{tikzpicture}[scale=.75,line width = 1.5\Lpt,vertex/.style={minimum size=50\Spt,draw,circle,thick},auto]

\tikzset{%
  % use: dblarw={color}{totalouterwidth}{outlineinsidewidth}
  dblarw/.style n args={3}{%,
    -latex,
    line width=#2,
    draw=black,  % this draw and color could
    color=black, % be an additional style arg
    opacity=1,
    % note: just color=#1 makes all black! 
    % fill=#1 makes insides (in an open! curve) filled too!
    % draw=#1,color=#1, seems to work, though
    postaction={
      draw=#1,
      color=#1,
      line width=#2-#3,
      shorten >=2*#3,
      %shorten <=#3,
    },
    dblarw/.default={gray}{5pt}{1pt},
    dblarw/.initial={gray}{5pt}{1pt},
  }
}

\node[vertex,fill=hlt2Fill] at (0,10) (w) {\Large $w$};
\node[vertex,fill=hlt2Fill] at (0,7) (x) {\Large $x$};
\node[vertex] at (0,4) (y) {\Large $y$};
\node[vertex] at (0,1) (z) {\Large $z$};
\node[vertex] at (5,5.5) (p) {\Large $p$};

\draw[loosely dotted, line width = 0.5mm] (0,5.5) ellipse (1.5cm and 6cm);

\draw[loosely dotted, line width = 0.5mm] (5,5.5) ellipse (1.5cm and 6cm);

\node[circle](g1) at (1.5,-.5) {};
\node[left=of g1] {\Large $G_1$};

\node[circle](g2) at (6.5,-.5) {};
\node[left=of g2] {\Large $G_2$};

\draw[ultra thick, loosely dashed] (w) to (x);
\draw[dblarw={black}{2pt}{.5pt}] (x) -- (y) node[midway] {\large $t=2$};
\draw[ultra thick,loosely dashed] (y) to (z);
\draw[ultra thick] (x) to (p);
\draw[ultra thick] (y) to (p);
\draw[dblarw={black}{2pt}{.5pt}] (w) -- (p) node[midway,sloped] {\large $t=1$};
\draw[dblarw={black}{2pt}{.5pt}] (p) -- (z) node[below, midway,sloped] {\large $t\geq2$};

%\node at (2,-1) {$P_4$ induced};
\end{tikzpicture}}
\caption{Illustration of the proof of  Lemma \ref{Join2ForbiddenInducedSG}.  The dotted line from $w$ to $x$ indicates an edge except for the case of $K_2 \du 2K_1$.  The dashed line from $y$ to $z$ is only an edge in the case of $P_4$.}
\label{Fig:Forbidden}
\end{figure}
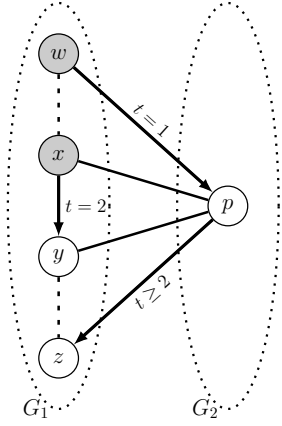

\begin{proof}
 Suppose that $G=G_1\vee G_2$ and at least one of the conditions
$ \PTpbar(G)=1$ or $\PTbar(G)=1$ holds.  The proof proceeds by assuming that $G$ has an induced subgraph $H$ isomorphic to one  of the three listed graphs, and showing that this implies a slow (standard or PSD) forcing set exists.  Figure \ref{Fig:Forbidden} illustrates the proof.
%There are six cases to be considered, depending on which of the three graphs is $H$ and depending on whether regular zero forcing or PSD zero forcing is to be considered.
%Two of these cases are already covered by the previous theorem, but the proof will cover all cases uniformly.

In all cases, the four vertices of $H$ must all lie within a single
side of the join, without loss of generality in $G_1$, because $H$ does not contain a $4$-cycle and has no vertex of degree greater than $2$.  In the case that $H=P_4$, label its vertices $(w, x, y, z)$ in order.
In the case that $H=P_3 \du K_1$, let the vertices of the $P_3$ be labeled
$(w, x, y)$ in order and let the vertex of the $K_1$ be labeled $z$.
In the case that $H=K_2 \du 2K_1$, let the vertices of the induced $K_2$
be labeled $x$ and $y$ and let the remaining vertices be labeled $w$ and $z$.
%\alert{Diagrams would be useful here, showing the three cases, with consistent labeling,
%and possibly illustrating $p$ in $G_2$ as well and the complement of $\{p, y, z\}$ colored blue.}
Let $p$ be any vertex of $G_2$.
and let  $B=V(G)\setminus \{p, y, z\}$. 

No blue vertex can force  $y$ or $z$ in the first round, because the set of white vertices
is connected and every blue vertex is either in $G_1$, and also adjacent to the white vertex $p$, or is
in $G_2$, and is adjacent to both white vertices $y$ and $z$. The blue vertex $w$ is adjacent
to $p$ and to no other white vertex, so $w$ forces $p$ in first round, then in the second round the vertex
$x$ forces $y$. The vertex $z$ can be forced by $p$ in
the third round  if it was not already forced in the second round.
Thus $B$ is a slow forcing set of $G$.
\end{proof}

 The next result extends the power of Lemma \ref{Join2ForbiddenInducedSG}.

\begin{theorem}{\rm \cite[Theorem 23]{alf20}}
\label{ForbiddenCharacterization}
Let $ H$ be a connected graph. Then $ H$ is $\{P_4,\textsf{\textup{paw}},\textsf{\textup{diamond}}\}$-free.
if and only if $H \cong K_n$ for some $n>0$ or $ H\cong K_{m,n}$ for some $m,n>0$.
\end{theorem}

\begin{corollary}
\label{CompleteOrCompleteBipartite}
If $G$ is a join of two graphs and at least one of the two conditions
$ \PTpbar(G)=1$ or $\PTbar(G)=1$ holds, then all of the components of $\overline{G}$ are complete graphs or complete bipartite graphs.
\end{corollary}

\begin{proof}
Let $G$ be a join of two graphs and suppose that at least one of
$ \PTpbar(G)=1$ or $\PTbar(G)=1$ holds.  By Lemma \ref{Join2ForbiddenInducedSG}, none of $P_4$, $P_3 \du K_1$,
or $K_2 \du 2K_1$ is an induced subgraph of $G$; therefore $\overline{G}$ is $\{P_4,\textsf{\textup{paw}},\textsf{\textup{diamond}}\}$-free.  Therefore, all of its components are as well, so by Theorem \ref{ForbiddenCharacterization}, they are complete graphs or complete bipartite graphs.
\end{proof}
 
 Now that we have preliminary tools that apply to both standard and PSD propagation time, we focus first on the PSD case. A set  $B\subseteq V(G)$ is \emph{independent} if $G[B]$ has no edges.   The \emph{independence number}  of $G$, denoted by $\alpha(G)$, is the maximum cardinality of an independent set.

\begin{lemma}
\label{Join2No3K1}
If $\PTbar_+(G)=1$ and $G$ is a join of two  graphs, then $\alpha(G)\le 2$  and therefore $\overline{G}$ has no induced subgraphs isomorphic to $K_3$.
\end{lemma}

\begin{proof}
Suppose $G=G_1\vee G_2$ and there are  three independent vertices  $u$, $v$, and $w$. These must all lie in one side of the join, without loss
of generality in $G_1$. Let $p$ be any vertex of $G_2$.  We show that the complement of the vertex set $\{p, u, v\}$ is a PSD forcing set that forces in time greater than one. No blue vertex
can force $u$ in the first round, because the set of white vertices is connected
and every blue vertex that is adjacent to $u$ is either in $G_1$, and is also adjacent to $p$, or is
in $G_2$ and is also adjacent to $v$. However, the blue vertex $w$ is adjacent to $p$ and no other white vertex, and so in the first round it forces $p$. Now the set of white vertices $\{u, v\}$ is
disconnected, and in the second round, $p$ can force both $u$ and $v$.
\end{proof}

%\bs{Lemma 5.4 implies, with Theorem 23 from \url{https://dx.doi.org/10.1016/j.laa.2019.09.012} (which is simple enough to prove directly if desired), that (if $G\not\cong K_n$) the components of $\overline{G}$ are either complete bipartite graphs or complete graphs.  The following lemma can be eliminated, since the preceding lemma eliminates the complete graph possibility.}

\begin{theorem}
\label{PSDJoin2ConjForm}
Let $\PTpbar(G)=1$ and let $G$ be a join of two  graphs. Then $G$ is a PSD fast join.
\end{theorem}

 \begin{proof}
 Since $G$ is the join of two graphs, $\Gc$ has at least two connected components.  By  Corollary \ref{CompleteOrCompleteBipartite} and Lemma \ref{Join2No3K1}, every component of $\Gc$ is a complete bipartite graph.  Consider  a component of $\overline{G}$ and let  $A$ and $B$ denote the partite sets.    If  $A$ (or $B$) is empty, then $B$ (or $A$) contains an isolated vertex of $\overline{G}$ and therefore a universal vertex of $G$, implying that $G\cong K_n$ (by  Corollary \ref{NoDomVtx}).  
\end{proof}

Next we consider standard forcing.
\begin{lemma}
\label{StdNoStars}
Let $G$ be a join of two  graphs $G_1$ and $G_2$, and $\overline{G_1}\cong K_{1,r}$ with $r\ge 2$. Then $\PTbar(G)\geq 2$.
\end{lemma}

\begin{proof}
%Suppose $G$ is a join of two  graphs $G_1$ and $G_2$, $k>1$, and and $\overline{G_1}\cong K_{1,r}$.  
Let $c$ be the center of the star $\overline{G_1}$, let $u$ and $v$ be leaves of $\overline{G_1}$,  let $w\in V(G_2)$, and define $B=V(G)\setminus\{c,v\}$.  
No vertex in $G_2$ can perform a force in the first round since such a vertex is  adjacent to both $v$ and $c$ in $G$. Since   $u$ is adjacent in $G$ to $v$ but not to $c$, $u$ can force $v$ in the first round.  However, since $N_G(c)=V(G_2)$, $c$ cannot be forced in the first round.  After round one,  $c$ is the only white vertex, so any of its neighbors (e.g., $w$) can then force $c$.
\end{proof}

\begin{theorem}
\label{StdJoin2Conjecture} 
If $\PTbar(G)=1$ and $G$ is a join of two  graphs, then $G$ is a standard fast join. 
\end{theorem}

\begin{proof}
Assume $\PTbar(G)=1$ and $G$ is a join of two  graphs, so $\Gc$ has at least two connected components.
If $G$ is not only a join but a complete graph, then $G$ is indeed a standard fast join.
Otherwise it suffices to show, for  any component $H$ of $\overline{G}$, that $H$ is either complete of order at least two, or
is complete bipartite with each partite set containing at least two vertices.

By  Theorem \ref{ForbiddenCharacterization}, $H$ is either complete or complete bipartite.
If $H$ is complete, it must be of order at least two so as to avoid the universal vertices that are forbidden by
 Corollary \ref{NoDomVtx}.
This leaves only the case where $H$ is complete bipartite.
The case where one partite set  has exactly one vertex and the other has more than one makes $H$
a star graph of the sort forbidden by  Lemma \ref{StdNoStars}.
Thus $A$ and $B$ {each} contain at least two vertices.  
\end{proof}

 The next two results allow us to resolve Conjectures \ref{c:psd} and \ref{c:std} for graphs $G$ with the associated forcing number at least $n-2$.

\begin{theorem}\label{compGZn-2}{\rm \cite[Theorem 5.1]{AIM}, \cite[Theorem 14]{BHL04}}  A graph $G$ of order $n$ has $\Z(G)\ge n-2$ if and only if the complement of $G$
is of the form 
\[(K_{s_1}\du \dots\du K_{s_t}\du K_{p_1,q_1}\du\cdots\du K_{p_k,q_k})\vee K_r\] for  appropriate
nonnegative integers $t, k,s_1,\dots,s_t,p_1,q_1,\dots,p_k,q_k,r$ with $p_i+q_i>0$ for all $i=1,\dots,k$. 
\end{theorem}

\begin{corollary}  Let $G$ be a connected graph of order $n\ge 2$.
If $\Z(G)\ge n-2$ and $\PTbar(G)=1$, then $G$ is a standard fast join.  If $\Zp(G)\ge n-2$ and $\PTpbar(G)=1$, then $G$ is a PSD fast join.
\end{corollary}
\bpf   Assume   $\Z(G)\ge n-2$ or $\Zp(G)\ge n-2$; the latter implies $\Z(G)\ge n-2$.  We apply Theorem \ref{compGZn-2} to show $G$ is a  join of two graphs: Observe first that since  $G$ is connected, $r=0$, i.e., $\Gc=K_{s_1}\du \dots\du K_{s_t}\du K_{p_1,q_1}\du\cdots\du K_{p_k,q_k}$.  If $t\ge 1$ and $k\ge 1$, or $t\ge 2$,  or $k\ge 2$ then $G$ is a join. If $k=0$ then $t\ge 2$ since $G$ is connected and $n\ge 2$. So suppose $t=0$ and $k=1$.  Note that  $\Gc=K_{p_1,q_1}$ with $p_1,q_1\ge 1$ would result in a disconnected graph $G$.  So $p_1\ge 2$ and $q_1=0$, and $G=K_{p_1}$ is a complete graph (which is a join since $n\ge 2$).  
  Since $G$ is a join, the statements  follow from  Theorems \ref{PSDJoin2ConjForm} and \ref{StdJoin2Conjecture}.
\epf

%\bigskip
 \section*{Acknowledgments}

 The authors thank the American Institute of Mathematics (AIM) and the National Science Foundation (NSF) for support of this research. 
The research of M. Hunnell and L. Hogben was also partially supported by NSF grant 2447261.

%==========================
%
% Section
%
%==========================

\end{document}